\documentclass[11pt]{amsart}
\usepackage{amscd}
\usepackage{amssymb}
\usepackage[matrix,arrow]{xy}
\usepackage{color}

\newtheorem{theorem}{Theorem}[section]
\newtheorem{proposition}[theorem]{Proposition}
\newtheorem{lemma}[theorem]{Lemma}

\newtheorem{corollary}[theorem]{Corollary}
\theoremstyle{remark}

\newtheorem{examples}[theorem]{Examples}
\theoremstyle{example}

\newcommand{\id}{\mathrm{id}}

\newcommand{\Aut}{\mathrm{Aut}}
\newcommand{\Ext}{\mathrm{Ext}}
\newcommand{\Hom}{\mathrm{Hom}}

\newcommand{\ot}{\otimes}

\newcommand{\KK}{\mathrm{KK}}
\newcommand{\KKX}{\mathrm{KK}_{X}}

\newcommand{\Tor}{\mathrm{Tor}}

\newcommand{\ho}{homomorphism}
\newcommand{\mo}{monomorphism}

\newcommand{\mos}{monomorphisms}

\newcommand{\CCC}{\mathbb{C}}

\newcommand{\ep}{\varepsilon}

\newcommand{\QQQ}{\mathbb{Q}}
\newcommand{\ZZZ}{\mathbb{Z}}

\newcommand{\OOO}{\mathcal{O}}

\newcommand{\fset}{\mathcal{F}}

\newcounter{rocount}

\begin{document}

\title{Fiberwise $KK$-equivalence of continuous fields of C*-algebras }
\thanks{The author was partially supported by
NSF grant \#DMS-0500693}
\maketitle

 \centerline{\normalsize
MARIUS DADARLAT} \vskip 4pt  \centerline{\it \footnotesize Purdue University}
\centerline{\it\footnotesize West Lafayette, IN, U.S.A.} \vskip 4pt

%


\begin{abstract}

Let $A$ and $B$ be separable nuclear continuous $C(X)$-algebras
over a finite dimensional compact metrizable space $X$.
It is shown that an element $\sigma$ of the parametrized Kasparov
 group $\KKX(A,B)$ is invertible if and only all if its fiberwise components
$\sigma_x\in KK(A(x),B(x))$ are invertible.
This criterion does not extend to infinite dimensional spaces since
there exist nontrivial unital separable continuous fields
over the Hilbert cube with all fibers isomorphic to the Cuntz algebra $\OOO_2$.
Several applications to continuous fields of Kirchberg algebras are given. It is
also shown that if each fiber of a
separable nuclear continuous $C(X)$-algebra $A$
over a finite dimensional locally compact space $X$  satisfies
the UCT, then $A$ satisfies the UCT.
    \end{abstract}


\section{Introduction}
Continuous C*-bundles arise naturally:
any separable C*-algebra $A$ with Hausdorff primitive space $X$
is isomorphic to the C*-algebra of continuous sections of a
 C*-bundle  over $X$  with fibers  the primitive quotients of $A$ \cite{Fell}, \cite{BK:bundles}.
A continuous C*-bundle (also called continuous field \cite{Fell} or continuous $C(X)$-algebra \cite{Kas:inv}) needs not be locally trivial at any point. In his work on the Novikov conjecture \cite{Kas:inv}, Kasparov has introduced
 parametrized KK-theory groups $\mathcal{R}KK(X;A,B)$  for $C(X)$-algebras $A$ and $B$.
These groups, for which we prefer the more compact notation $KK_X(A,B)$, admit a natural product structure $KK_X(A,B)\times KK_X(B,C)\to KK_X(A,C)$. The invertible elements in $KK_X(A,B)$
are denoted by $KK_X(A,B)^{-1}$. If $KK_X(A,B)^{-1}\neq \emptyset$
we say that $A$ is $KK_X$-equivalent to $B$. The work of Kirchberg
on the classification of purely infinite C*-algebras raises the question
of determining when two $C(X)$-algebras are $KK_X$-equivalent. Indeed, by \cite{Kir:Michael}, if $A$ and $B$ are two unital separable nuclear
C*-algebras with Hausdorff spectrum $X$, then $A\otimes \OOO_\infty \otimes \mathcal{K}\cong B\otimes \OOO_\infty \otimes \mathcal{K}$ if and only if $A$ is $KK_X$-equivalent to $B$.
A closely related problem is to characterize the invertible
elements of $\KKX(A,B)$. We give a criterion for  $\KKX$-invertibility
  for spaces $X$
of finite covering dimension.

\begin{theorem}\label{Theorem:general-KK}
   Let $A$ and $B$ be separable nuclear continuous $C(X)$-algebras
over a finite dimensional compact metrizable space $X$.
   If $\sigma\in \KKX(A,B)$, then  $\sigma\in \KKX(A,B)^{-1}$ if and only if
   $\sigma_x\in
   KK(A(x),B(x))^{-1}$ for all $x\in X$.
\end{theorem}
Consequently, if the fibers of $A$ and $B$ satisfy the Universal Coefficient
Theorem for the Kasparov groups (abbreviated UCT, \cite{RosSho:UCT}) and if there is a $C(X)$-linear morphism from $A$ to $B$, or just an element of
$\KK_X(A,B)$, such that all the  maps $K_*(A(x))\to K_*(B(x))$
are bijective, then $K_*(A)\cong K_*(B)$. This opens the way for the use of homological methods for the computation of the K-theory groups
of continuous fields.
Let us note that Mayer-Vietoris type arguments are
not directly applicable due to lack of local triviality.
The assumption on the finite dimensionality of $X$ is essential:
\begin{examples}\label{example:nontriv-O-two}
 There is a family $(E_P)_P$ which has
 the power of continuum and which consists of
 mutually nonisomorphic unital separable
continuous $C(Z)$-algebras over the Hilbert cube $Z$ with all fibers isomorphic to
the Cuntz algebra $\OOO_2$.
\end{examples}
The $K_0$-groups of $E_P$ are nonzero even though
$Z$ is contractible
and  $KK(\OOO_2,\OOO_2)=0$.
The family $(E_P)_P$ is easily constructed starting from
an example of a continuous $C(Y)$-algebra over $Y=\prod_{n=1}^\infty S^2$
with fibers isomorphic to the CAR algebra but which does not
absorb the CAR algebra, exhibited in \cite{HirshbergRordamWinter:absorb-ssa}.

By specializing Theorem~\ref{Theorem:general-KK} to the case when $A$ is a trivial $C(X)$-algebra,  we obtain, based on \cite{Kir:Michael}, an explicit necessary and sufficient (Fell type) K-theory condition for triviality
of continuous fields of arbitrary Kirchberg algebras: Corollary~\ref{Cor:Fell-cond}.
 In particular we no longer require KK-semiprojectivity of the fibers as in our earlier approach \cite{Dad:bundles-fdspaces}.

Let us say that a unital C*-algebra $D$ has
 the \emph{automatic  triviality
property}   if any separable unital
continuous $C(X)$-algebra over
a finite dimensional
 compact
 metrizable  space $X$ all of whose fibers are isomorphic
 to $D$ is isomorphic to $C(X)\ot D$.
We proved in \cite{Dad:bundles-fdspaces} that the Cuntz algebras
$\OOO_2$ and  $\OOO_\infty$ are the only
Kirchberg algebras with the automatic triviality
property among those Kirchberg algebras satisfying the UCT and having
finitely generated K-theory.
  By combining
Corollary~\ref{Cor:Fell-cond} with
our homotopy calculations from \cite{Dad:homotopy-aut} and the recent absorption
result of Hirshberg, R{\o}rdam and Winter \cite{HirshbergRordamWinter:absorb-ssa}, we are now able to drop the finite generation condition.
\begin{theorem}\label{Theorem:ssa-Kirchberg-algebras}
  A unital Kirchberg algebra $D$ satisfying the UCT has the automatic  triviality
property if and only $D$ is isomorphic to either $\OOO_2$,  $\OOO_\infty$, or
$\OOO_\infty\ot U$ where $U$ is a unital uniformly hyperfinite algebra of infinite type.
\end{theorem}
The condition that $X$ is finite dimensional cannot be dropped
as shown by Example~\ref{example:nontriv-O-two}.
Let us recall that a Kirchberg algebra is a purely infinite simple nuclear separable C*-algebra
\cite{Ror:ency} and that a separable C*-algebra $A$
satisfies the UCT  if and only if $A$ is KK-equivalent to a commutative C*-algebra \cite{RosSho:UCT}. The class of C*-algebras satisfying the UCT is surprisingly large;
it contains the C*-algebras of
 locally compact second countable amenable groupoids  \cite{Tu:BC}.
\cite[Thm.~4.3]{HirshbergRordamWinter:absorb-ssa}.
The list given by Theorem~\ref{Theorem:ssa-Kirchberg-algebras}
coincides with the list of all strongly selfabsorbing Kirchberg algebras
satisfying the UCT exhibited in \cite{WinterToms:ssa}.
In a forthcoming joint paper with  Winter \cite{Dadarlat-Winter:KK-of-ssa}, we
show that any $K_1$-injective separable strongly selfabsorbing
C*-algebra has the automatic triviality
property.

In the last part of the paper we prove a new
permanence property for the class of nuclear C*-algebras which satisfy the UCT:
\begin{theorem}\label{intro-uct}
    A separable nuclear continuous $C(X)$-algebra
   over a finite dimensional locally compact  space
   satisfies the UCT
 if all its fibers satisfy the UCT.
\end{theorem}


The  author is grateful to   {\'E}.~Blanchard for useful discussions and comments. He is also indebted to S. Echterhoff for a conversation on triviality criteria for continuous fields which highlighted the role of continuous fields with
KK-contractible fibers. This inspired the author to prove
Theorem~\ref{Theorem:general-KK}.

\section{$C(X)$-algebras and $KK_{X}$-equivalence}\label{KK_{X}-equivalence}

 Kirchberg
has shown that any nuclear separable C*-algebra is equivalent in KK-theory to a
Kirchberg algebra \cite[Prop.~8.4.5]{Ror:ency}. We extend his
result in the context of continuous $C(X)$-algebras and $KK_{X}$-theory (see
Theorem~\ref{C(X)-K-version}). The space $X$ is assumed to be compact and metrizable
throughout this section.
\begin{lemma}[{\cite[Prop.~3.2]{Blanchard:Hopf}}]\label{Ph-K-unitization}
    If $A$ is a continuous $C(X)$-algebra, then
    there is a split
    short exact sequence of continuous $C(X)$-algebras
\begin{equation}\label{Blanachard-unitization}
    \xymatrix{
 {0}\ar[r]& {A}\ar[r]
                & {A^+} \ar@<0.5ex>[r]&C(X)\ar@<0.5ex>[l]^{\alpha}\ar[r]&0
}
\end{equation}

where $A^+$ is unital, $\alpha$ is $C(X)$-linear and $\alpha(1)=1$.
\end{lemma}
Consider the category of separable $C(X)$-algebras where the morphisms from $A$
to $B$ are the elements of $KK_{X}(A,B)$ with composition given by the
Kasparov product. The isomorphisms in this category are the
$\KKX$-invertible
elements denoted by $KK_{X}(A,B)^{-1}.$ Two $C(X)$-algebras are
$KK_{X}$-equivalent if they are isomorphic objects in this category. In the
sequel we shall use twice the following elementary observation (valid in any
category). If composition with $\gamma \in KK_{X}(A,B)$ induces a bijection
$\gamma^*:KK_{X}(B,C)\to KK_{X}(A,C)$ (or $\gamma_*:KK_{X}(C,A)\to KK_{X}(C,B)$) for  $C=A$ and $C=B$, then $\gamma\in
KK_{X}(A,B)^{-1}$.
\begin{lemma}\label{O2}
Let $A$ be a  separable nuclear continuous $C(X)$-algebra.
      Then there exist a  separable nuclear unital continuous $C(X)$-algebra
      $A^\flat$  and
     $C(X)$-linear \mos\ $\alpha:C(X)\otimes \OOO_2 \to
     A^\flat$  and $\jmath :A \to A^\flat$ such that $\alpha$ is
      unital
     and $\KKX(\jmath)\in \KKX(A,A^\flat)^{-1}$.
\end{lemma}
\begin{proof} Let $p\in \OOO_\infty$ be a non-zero projection with
$[p]=0$ in $K_0(\OOO_\infty)$. Then there is a unital $*$-\ho\ $\OOO_2\to
p\OOO_\infty p$ which induces a $C(X)$-linear unital $*$-\mo\
$\mu:C(X)\ot\OOO_2\to C(X)\ot p\OOO_\infty p$. We tensor the exact sequence
\eqref{Blanachard-unitization} by $p\OOO_\infty p$ and then take the pullback
by $\mu$. This gives a split exact sequence of unital $C(X)$-algebras:
\[
\xymatrix{
 {0}\ar[r]& {A\ot p\OOO_\infty p}\ar[r]\ar@{=}[d]
                & {A^+\ot p\OOO_\infty p} \ar@<0.5ex>[r]&C(X)\ot
                 p\OOO_\infty p\ar@<0.5ex>[l]^{\alpha}\ar[r]&0\\
{0}\ar[r]& {A\ot p\OOO_\infty p}\ar[r]^{j}
                & {A^\flat}\ar[u] \ar@<0.5ex>[r]&
                C(X)\ot \OOO_2\ar[u]_{\mu}\ar@<0.5ex>[l]^{\alpha}\ar[r]&0
 }\]
The map $A^\flat \to {A^+\ot p\OOO_\infty p}$ is a unital $C(X)$-linear
$*$-monomorphism, so that $A^\flat$ is a continuous $C(X)$-algebra. It is
nuclear being an extension of nuclear C*-algebras.  By \cite[Thm.
5.4]{Bauval:KKX} for any separable nuclear continuous $C(X)$-algebra $B$ there is an exact
sequence of groups
%

\[ \xymatrix{
{0\to \KKX(B,A\ot p\OOO_\infty p)}\ar[r]^-{j_*}
                & {\KKX(B,A^\flat)} \ar[r]&
                {\KKX(B,C(X)\otimes \OOO_2)\to 0}.
 }\]

 $\KKX(B,C(X)\ot \OOO_2)=0$ since the class of the identity
  map of $C(X)\ot \OOO_2$ vanishes in
 $\KKX$. Therefore $j_*$ is bijective and so
  $\KKX(\jmath)\in \KKX( A\otimes p\OOO_\infty p,A^\flat)^{-1}$.
We conclude
the proof by observing that map  $A \to A\otimes p\OOO_\infty p$, $a\mapsto a
\otimes p$, induces a $\KKX$-equivalence.
 \end{proof}
 Let $(A_i,\varphi_i)$ be an inductive system
 of separable nuclear unital continuous $C(X)$-algebras
   with unital (fiberwise) injective connecting maps. Let $A=\varinjlim\, (A_i,\varphi_i)$ be
   the inductive limit C*-algebra and let
   $\varphi_{i,\infty}:A_i\to A$ be the induced inclusion map.

 \begin{lemma}\label{Lemma:admisible-system} $A$ is a unital nuclear continuous $C(X)$-algebra
 and there is a sequence $(\eta_i:A\to A_{n(i)})_i$ of unital completely positive
   $C(X)$-linear maps such that $(\varphi_{n(i),\infty}\circ \eta_i)_i$
   converges to $\id_A$ in the point norm topology.
 \end{lemma}
 \begin{proof} The map $C(X)\to Z(A)$ is induced by the maps
 $C(X)\to Z(A_i)$, so that $A$ clearly becomes a $C(X)$-algebra.
 The continuity of the map $x\mapsto \|a(x)\|$ for $a\in A$
 is verified by approximating $a\in A$ by some $a_i\in A_i$ and using
 the fiberwise injectivity of $\varphi_{i,\infty}$.
 Since $A$ is a nuclear C*-algebra, it follows by \cite[Thm.~7.2]{Bauval:KKX}
 that $A$ is $C(X)$-nuclear. This means that there are sequences of
 unital $C(X)$-linear completely positive  maps $\alpha_i:A\to C(X)\ot M_{k(i)}$
 and $\beta_i:C(X)\ot M_{k(i)}\to A$ such that $\beta_i\circ\alpha_i$
 converges to $\id_A$ in the point-norm topology.
 After perturbing the restriction of $\beta_i$ to $1_{C(X)}\otimes M_{k(i)}$
 to a unital and completely positive map  $\beta'_i:1_{C(X)}\otimes M_{k(i)}\to A_{n(i)}$
 and extending $\beta'_i$ to a $C(X)$-linear map
  $\beta''_i:C(X)\ot M_{k(i)}\to A_{n(i)}$, we may assume that
$\beta_i$ factorizes as $\beta_i=\varphi_{n(i),\infty}\circ \beta''_i$.
Then $\eta_i=\beta''_i\circ \alpha_i$ satisfies the conclusion of the lemma.
 \end{proof}

\begin{proposition}\label{cor-milnor}
 Let $(A_i,\varphi_i)$ be an inductive system of separable nuclear unital
 continuous
  $C(X)$-algebras
   with unital injective connecting maps.
    If $\varphi_i\in \KKX(A_i,A_{i+1})^{-1}$
    for all $i$, and $\Phi:A_1 \to \varinjlim
   (A_i,\varphi_i)=A$ is the induced map,
   then $\Phi\in \KKX(A_1,A)^{-1}$.
\end{proposition}
\begin{proof}
We use  Milnor's $\varprojlim^1$-exact sequence in $\KKX$-theory applied to
the inductive system $(A_i,\varphi_i)$. The proof of $\sigma$-additivity of
$KK(A,B)$ in the first variable given in \cite[Thm.~2.9]{Kas:inv} applies with
essentially no changes to show the corresponding property for  $\KKX(A,B)$.
Lemma~\ref{Lemma:admisible-system} verifies the assumptions of
\cite[Lemma~2.7]{Meyer-Nest:BC}. Thus the system $(A_i,\varphi_i)$ is admissible
in the sense of \cite[Def.~2.5]{Meyer-Nest:BC} and hence by \cite[Lemma~2.4 and
Prop.~2.6]{Meyer-Nest:BC} we have an exact sequence:
\[
\xymatrix{
 {0\to\varprojlim^1 \KKX^{1}(A_i,B)}\,\ar[r]
                & {\KKX(A,B)}\ar[r]&
                {\varprojlim\, \KKX(A_i,B)\to 0}
               }
\]
(One can also give a direct proof of the exact sequence which is essentially
identical to the proof of the corresponding sequence in KK-theory. One argues as
in \cite{RosSho:UCT} using the exact sequences from \cite{Bauval:KKX}.  The maps
from Lemma~\ref{Lemma:admisible-system} are needed to verify that the mapping
telescope extension of $A$ is semisplit in the category of $C(X)$-algebras.)

 Since $\varprojlim^1 \big(G_i,\lambda_i\big)=0$ and $ G_1\cong\varprojlim
 \big(G_i,\lambda_i)$
 for any sequence of abelian groups $(G_i)_{i=1}^\infty$ and
 group isomorphisms $\lambda_i:G_i\to G_{i+1}$,
 the $\varprojlim^1$-exact sequence shows
  that for any separable continuous $C(X)$-algebra $B$,
  the map
 $\KKX(A,B)\rightarrow \KKX(A_1,B)$ induced by $\Phi$ is
 bijective.
 Therefore $\KKX(\Phi)\in \KKX(A_1,A)^{-1}.$
\end{proof}

 We need the following $C(X)$-equivariant construction which
 parallels a construction of Kirchberg as presented in
 \cite{Ror:ency}. A similar deformation technique has appeared
 in \cite{Dad:rfd1}.
 \begin{theorem}\label{C(X)-K-version}
     Let $A$ be a  separable nuclear
     continuous $C(X)$-algebra. Then  there exist  a separable nuclear
     continuous unital $C(X)$-algebra $A^\sharp$ whose fibers are Kirchberg
     C*-algebras and a $C(X)$-linear $*$-\mo\ $\Phi:A \to A^\sharp$
     such that $\Phi$ is a $\KKX$-equivalence.
     For any $x\in X$ the map $\Phi_x:A(x)\rightarrow A^\sharp(x)$ is a
     $KK$-equivalence.
     \end{theorem}
\begin{proof} By Proposition~\ref{O2} we may assume that $A$ is unital and that
there is a unital $C(X)$-linear $*$-\mo\ $\alpha:C(X)\otimes \OOO_2 \to A$. By
\cite[Thm. 2.5]{Blanchard:subtriviality} there is a unital $C(X)$-linear
$*$-monomorphism \mbox{$\beta:A \to C(X)\ot\OOO_2$}. Let $s_1,s_2$ be the
images in $A$  of the canonical generators $v_1,v_2$ of  $\OOO_2\subset
C(X)\ot\OOO_2$ under the map $\alpha$. Set $\theta=\alpha\circ \beta:A \to A$
and define  $\varphi:A \to A$ by $\varphi(a)=s_1\, a\,
s_1^*+s_2\,\theta(a)\,s_2^*$. The unital $*$-\ho\ $\varphi_x:A(x)\to A(x)$
induced by $\varphi$ satisfies $\varphi_x\pi_x=\pi_x\varphi$ and
$\varphi_x(b)=s_1(x)\, b\, s_1(x)^*+s_2(x)\,\theta_x(b)\,s_2(x)^*$.
 Moreover $\theta_x$ factors through $\OOO_2$ since
  $\theta_x=\alpha_x\circ \beta_x$.
 Let $A^\sharp$ be the continuous $C(X)$-algebra obtained as the
 limit of the inductive system
\[\xymatrix{
{A}\ar[r]^{\varphi}&{A}\ar[r]^{\varphi}& {A}\ar[r]^{\varphi}& {\cdots}
 }\]
 and let $\Phi:A \to A^\sharp$ be the induced map.
 The commutative diagram
 \[\xymatrix{
{A}\ar[d]^{\pi_x}\ar[r]^\varphi& {A}\ar[d]^{\pi_x}\ar[r]^\varphi &
{A}\ar[d]^{\pi_x}
\ar[r]^\varphi&{\cdots}\ar[r]&{A^\sharp}\ar[d]^{\pi_x}\\
{A(x)}\ar[r]^{\varphi_x}& {A(x)}\ar[r]^{\varphi_x}&
{A(x)}\ar[r]^{\varphi_x}&{\cdots}\ar[r]&{A^\sharp(x)}
 }\]
 shows that the fiber   $A^\sharp(x)$ of $A$ is isomorphic to $\varinjlim
 (A(x),\varphi_x)$.
By the proof of \cite[Prop.~8.4.5]{Ror:ency} $A^\sharp(x)$ is a unital
Kirchberg algebra. It remains to prove that the  map $\Phi:A \to A^\sharp$
induces a $\KKX$-equivalence. By Proposition~\ref{cor-milnor} it suffices
to verify that $\KKX(\varphi)=\KKX(\id_A)$. This follows from the
equation $\varphi(a)=s_1\, a \,s_1^*+s_2\,\theta(a)\,s_2^*$,   since $\theta$
factors through $C(X)\otimes \OOO_2$ and hence $\KKX(\theta)=0$.
\end{proof}
\begin{theorem}\label{Lemma:contactible-fibers}
   Let $X$ be a compact metrizable finite dimensional space.
   Let $A$  be a separable nuclear continuous $C(X)$-algebra.
   If all the fibers of $A$ are $\KK$-contractible,
   then $A$ is $\KKX$-contractible.
\end{theorem}
\begin{proof} By Theorem~\ref{C(X)-K-version}, $A$ is $\KKX$-equivalent
to  a separable nuclear unital continuous $C(X)$-algebra $A^\sharp$ whose fibers
are $\KK$-contractible Kirchberg algebras. Therefore $A^\sharp(x)\cong \OOO_2$
for all $x$ \cite{Ror:ency}.
By
\cite[Thm.~1.1]{Dad:bundles-fdspaces} $A^\sharp$ is
isomorphic to $C(X)\otimes\OOO_2$ and hence is $\KKX$-contractible.
Alternately one can argue that $A^\sharp\cong A^\sharp\otimes \OOO_2$
by \cite{HirshbergRordamWinter:absorb-ssa} and hence that
$\KKX(A^\sharp,A^\sharp)=0$.
\end{proof}
\emph{Proof of Theorem~\ref{Theorem:general-KK}}
\begin{proof} By Theorem~\ref{C(X)-K-version} we may assume that both $A$ and
$B$ are stable continuous  $C(X)$-algebras which absorb
$\OOO_\infty$ tensorially and whose fibers are Kirchberg algebras. 
By \cite[Hauptsatz 4.2]{Kir:Michael}, for any given $\sigma\in \KKX(A,B)$,
 there is a $C(X)$-linear $*$-\ho\
$\varphi:A\to B$ such that $\KKX(\varphi)=\sigma$. The mapping cone of
$\varphi$,
\[C_\varphi=\{(a,f)\in A\oplus C_0[0,1)\ot B\,:\,  f(0)=\varphi(a)\} \]
is a separable nuclear continuous $C(X)$-algebra with fibers
$C_{\varphi}(x)\cong C_{\varphi_x}$, $x\in X$. Since each $\varphi_x$ is a
KK-equivalence, it follows from the Puppe sequence in KK-theory
\cite{Bla:k-theory} that $C_{\varphi_x}$ is KK-contractible for each $x\in X$.
Then $C_\varphi$ is $\KKX$-contractible by
Theorem~\ref{Lemma:contactible-fibers}. Using   the Puppe exact sequence for
separable nuclear continuous $C(X)$-algebras (see \cite{Bauval:KKX})
\[ \xymatrix{
 \KKX(C,C_\varphi)\ar[r]
                & \KKX(C,A)
                \ar[r]^-{\varphi_*}&\KKX(C,B)\ar[r]&{\KKX^1(C,C_\varphi)}
 }\]
 we see that
$\varphi_*:\KKX(C,A)\to \KKX(C,B)$
is bijective for all separable nuclear
continuous $C(X)$-algebras $C$ and hence $\sigma$ is a $\KKX$-equivalence.
\end{proof}
A  remarkable isomorphism result for  separable
 nuclear strongly purely infinite
stable C*-algebras
 was announced (with an outline of the proof) by
Kirchberg in \cite{Kir:Michael}: two such C*-algebras $A$ and $B$ with the same
primitive spectrum $X$ are isomorphic if and only if they are
$\KKX$-equivalent.  In conjunction with Theorem~\ref{Theorem:general-KK} we derive the following.
\begin{theorem}\label{Theorem:general-KK-for Kirchberg algebras}
   Let $X$ be a compact metrizable finite dimensional space.
   Let $A$ and $B$ be separable continuous
   $C(X)$-algebras all of whose  fibers are Kirchberg algebras.
   Suppose that there is $\sigma\in \KKX(A,B)$ such that $\sigma_x\in
  KK(A(x),B(x))^{-1}$ for all $x\in X$. Then there is an isomorphism
  of $C(X)$-algebras $\varphi:A\ot \mathcal{K} \to B\ot \mathcal{K}$
  such that $\KKX(\varphi)=\sigma$. Moreover if $A$ and $B$ are unital
  and if $K_0(\sigma)[1_A]=[1_B]$, then $A\cong B$.
\end{theorem}
\begin{proof}
 Since $X$ is finite dimensional, $A\ot \mathcal{K}\ot \OOO_\infty
\cong A\ot \mathcal{K}$ and $B\ot \mathcal{K}\ot \OOO_\infty \cong B\ot
\mathcal{K}$ by \cite[Cor.~5.11]{BK}. In the unital case, one also has
 $A\ot \OOO_\infty\cong
A $ and $B\ot \OOO_\infty\cong B$ by \cite[Cor.~5.11]{BK}, \cite[Thm.
4.23]{KR1} and \cite[Thm. 8.6]{KR2} as explained for example in
\cite[Lemma~3.4]{DadPas:fields-over-zero-dim}. By
Theorem~\ref{Theorem:general-KK} $\sigma$ is a $\KKX$-equivalence. This enables
us to apply Kirchberg's  result \cite[Folgerung 4.3]{Kir:Michael} to
obtain an isomorphism of $C(X)$-algebras $\varphi:A\ot \mathcal{K} \to B\ot
\mathcal{K}$ such that $\KKX(\varphi)=\sigma$. In the unital case, since both
$\varphi(1_A\ot e_{11})$ and $1_B\ot e_{11}$ are full and properly infinite
projections in $B\ot \mathcal{K}$, the condition $\varphi_*[1_A]=[1_B]$ will
allows us to arrange that $\varphi(1_A\ot e_{11})=1_B\ot e_{11}$ after
conjugating $\varphi$ by a suitable unitary in $M(B\ot \mathcal{K})$ (see
~\cite{Cuntz:KOn}) and hence conclude that $A\cong B$.
\end{proof}
\begin{corollary}\label{Cor:Fell-cond}
   Let $X$ be a compact metrizable finite dimensional space.
   Let $A$  be separable continuous
   unital $C(X)$-algebra all of whose  fibers are Kirchberg algebras and let $D$ be a
   unital Kirchberg algebra.
   Suppose that there is $\sigma\in KK(D,A)$ such that $\sigma_x\in
  KK(D,A(x))^{-1}$ for all $x\in X$ and $K_0(\sigma)[1_A]=[1_B]$. Then  $A\cong C(X)\ot D$.
\end{corollary}
\begin{proof} This follows from the previous theorem since
$\KKX(C(X)\ot D,A)\cong KK(D,A)$.\end{proof}
\section{Nontrivial $\OOO_2$-bundles}\label{NT-bundles}
Let $Z$ denote the Hilbert cube.
Let $K$ be the Cantor set. If $G$ is a discrete group  let
$C(K,G)$ denote the continuous functions from $K$ to $G$.
For any two  countable, abelian, torsion groups $G_0$ and $G_1$,
we exhibit a unital separable continuous $C(Z)$-algebra $E$ with all fibers isomorphic to $\OOO_2$ such that $K_i(E)=C(K,G_i)$, $i=0,1$.

A crucial ingredient of the construction is an example from \cite{HirshbergRordamWinter:absorb-ssa} which we now recall
(with a minor variation).
Let $e$ be the unit of $C(S^2)$ and let $f\in M_2(C(S^2))$ be
the Bott projection. For each $n\geq 1$ we let $e_n=e$ and
$f_n=(e\oplus ... \oplus e) \oplus f$ ($n-1$ copies of $e$) be realized as orthogonal projections in $M_{n+2}(C(S^2))$.
Set
\[B_n=(e_n+f_n)M_{n+2}(C(S^2))(e_n+f_n), \quad A=\bigotimes_{n=1}^\infty B_n.\]

Let $U$ be the universal UHF  algebra with $K_0(U)=\QQQ$.
Let $Y=\prod_{n=1}^\infty S^2$. Arguing as in \cite{HirshbergRordamWinter:absorb-ssa} one shows that $A$ is a continuous
$C(Y)$-algebra with all fibers isomorphic to $U$.
Since each $B_n$ is Morita equivalent to $C(S^2)$,
  $K_0(B_n)$
 is freely generated as a $\ZZZ$-module by the classes of $[e_n]$
and $[f_n]$. It is clear that $K_1(A)=0$. Using the K\"unneth formula one shows
that $K_0(A)$ is isomorphic to the limit of the inductive system of groups
\[
\ZZZ^2 \to \ZZZ^4 \to \cdots\to \ZZZ^{2^n}\to \ZZZ^{2^{n+1}}\to \cdots
 \]
 where the $n$-th connecting morphism maps $x\in \ZZZ^{2^n}$ to  $(x,x)\in \ZZZ^{2^n}\oplus \ZZZ^{2^n} \cong \ZZZ^{2^{n+1}}.$
Consequently $K_0(A)$ is isomorphic to $C(K,\ZZZ)$.
Let $D$ be a unital Kirchberg C*-algebra such that
$K_i(D)=G_i$, $i=0,1$ and $D$ satisfies the UCT.
Then $F=A\otimes D$ is a unital separable continuous $C(Y)$-algebra
whose fibers are isomorphic to $A(x)\otimes D\cong U\otimes D.$
By the K\"unneth formula $K_i(U\otimes D)\cong \QQQ\otimes G_i=0$ since both $G_0$ and $G_1$ are torsion groups.
It follows that all the fibers $F(x)$ of $F$ are isomorphic to $\OOO_2$
 by the Kirchberg-Phillips classification theorem \cite{Ror:ency}.
On the other hand $K_i(F)=K_0(A)\otimes K_i(D) = C(K,\ZZZ)\otimes G_i\cong C(K,G_i)$, $i=0,1$.

By Blanchard's embedding theorem \cite{Blanchard:subtriviality}, there is a unital $C(Y)$-linear
monomorphism $\eta:F\to C(Y)\otimes \OOO_2$.
Let us regard $Y$ as a compact subset of the Hilbert cube $Z$.
Define
\[E=\{f\in C(Z,\OOO_2)\,:\, f|_Y \in \eta(F)\}.\]
Then $E$ is a separable unital continuous $C(Z)$-algebra with
all fibers isomorphic to $\OOO_2$.
Using the  exact sequence:
\[ \xymatrix{
0\ar[r]& C_0(Z\setminus Y, \OOO_2)\ar[r]& E
                \ar[r]&F \ar[r]& 0
} \]
we see that $K_i(E)\cong K_i(F)\cong C(K,G_i)$, $i=0,1$.
In particular $E$ is not isomorphic to $C(Z)\otimes \OOO_2$ if
$G_0\neq 0$ or $G_1\neq 0$, for example if $D=\OOO_n$, $2<n<\infty$.
For  a nonempty set $P$ of prime numbers, let $\ZZZ(P)$ be the subgroup
of $\mathbb{T}$ consisting of all elements whose orders have all prime
factors in $P$. One verifies that $C(K,\ZZZ(P))$ is not isomorphic to
$C(K,\ZZZ(P'))$ if $P\neq P'$. Indeed if $p\in P\setminus P'$, then
$C(K,\ZZZ(P'))$ does not have elements of order $p$, unlike $C(K,\ZZZ(P))$.
Consequently, the  family $(E_P)_P$ obtained by choosing $D$ with $K_0(D)=\ZZZ(P)$
consists of mutually non-isomorphic $C(Z)$-algebras and has the power of continuum.

\section{Automatic triviality}

\begin{proposition}\label{Tor}
Let $D$ be a unital Kirchberg algebra satisfying the UCT.
Suppose that $[X,\Aut(D)^0]$ reduces to singleton for any path connected
compact metrizable space $X$.
Then $D$ is isomorphic to either $\OOO_2$,  $\OOO_\infty$, or  $\OOO_\infty\ot U$
where $U$ is a unital UHF algebra of infinite type.
\end{proposition}

\begin{proof} We show that $D$ has the same K-theory
groups as one of the  listed  C*-algebras.
Let $C_{\nu}$ denote the mapping cone C*-algebra
of the unital inclusion $\nu:\mathbb{C} \to D$.
By \cite[Thm.~5.9]{Dad:homotopy-aut}
 there is a bijection
$[X,\Aut(D)^0]\to KK(C_{\nu},SC(X,x_0)\ot D)$
 for some (any) point $x_0\in X$.
Since the K-theory groups of $C(X,x_0)$ can be arbitrary countable
abelian groups it follows that $KK(C_{\nu},A\ot D)=0$ for all separable
C*-algebras $A$ satisfying the UCT.
Using the Puppe sequence in KK-theory
\cite{Bla:k-theory} we see that the restriction map
$\nu^*:KK(D,A\ot D)\to KK(\CCC,A\ot D)$ is bijective for all separable
C*-algebras $A$ satisfying the UCT. By the UCT (applied for  $A$ and its suspension), this implies that
(i) $K_1(D)=0$, (ii) $\Ext(K_0(D),K_0(A\ot D))=0$ and (iii) the  map
$\nu^*:\Hom(K_0(D),K_0(A\ot D))\to \Hom(K_0(\CCC), K_0(A\ot D))$ is bijective for all separable
C*-algebras $A$ satisfying the UCT.

First we are going to show that
 $G=K_0(D)$ is torsion free. We shall use the observation
that if $M,N$ are abelian groups and $M'$ is a subgroup of $M$ and $N'$ is a quotient of $N$, then $\Ext(M',N')$ is a quotient of $\Ext(M,N)$ \cite{Fuc:inf}.
Fix a prime $p$ and set $G[p]=\{x\in G:px=0\}$.
Using (ii) for $A=\OOO_{p+1}$, we obtain that $\Ext(G,G[p])=0$ since
$G[p]=\mathrm{Tor}(G,\ZZZ/p)$ is a quotient of $K_0(\OOO_{p+1}\ot D)$
and hence $\Ext(G,G[p])$ is a quotient of $\Ext(K_0(D),K_0(\OOO_{p+1}\ot D))=0$. Assuming that $G[p]\neq 0$ we find a subgroup of $G$ isomorphic
to $\ZZZ/p$ and hence $\Ext(\ZZZ/p,G[p])$ vanishes since it is a quotient of
$\Ext(G,G[p])=0$. On the other hand $\Ext(\ZZZ/p,G[p])$ is isomorphic
to $G[p]/pG[p]=G[p]$ and hence $G[p]=0$. This contradiction shows that
$G$ is torsion free.

Let $e$ denote the class of $[1_D]\in K_0(D)=G$.
If $e=0$, then $G=0$ by applying (iii) for $A=\CCC$.
For the rest of the proof we may assume that $e\neq 0$ and hence
that $\ZZZ\cong \ZZZ e\subset G$.

Let $A$ be a commutative C*-algebra such that $K_0(A)=G/\ZZZ$
\cite{Schochet:uct&milnor}. Since $G$ is torsion is free, by the K\"unneth
formula $K_0(A\ot D)\cong K_0(D)\ot K_0(A)\cong G\ot G/\ZZZ$. The exact sequence
$0\to \ZZZ\to G\to G/\ZZZ\to 0$ induces an exact sequence

$0\to \Hom(G/\ZZZ, K_0(A\ot D))\to \Hom(G, K_0(A\ot D))\to\Hom(\ZZZ, K_0(A\ot D))$.
From (iii) we obtain that $\Hom(G/\ZZZ, K_0(A\ot D))=0$ and hence
$\Hom(G/\ZZZ,G\ot G/\ZZZ)=0$.
Therefore for each $y\in G$, the morphism $G/\ZZZ\to G\ot G/\ZZZ$, $x\mapsto y\ot x$ is the zero map and hence $G\ot G/\ZZZ=0$. In particular $G/\ZZZ$ must be a torsion group since $G$ contains a copy of $\ZZZ$.
We also have that $\Tor(G/\ZZZ,G)=0$ since $G$ is torsion free.
The exact sequence $0\to \ZZZ\to G\to G/\ZZZ\to 0$ induces
an exact sequence
$\Tor(G/\ZZZ,G)\to \ZZZ\ot G \to G\ot G \to G/\ZZZ\ot G$.
 Therefore the map
$\theta:G \to G\ot G$, $\theta(x)=x\ot e$ is an isomorphism of groups. We also
know that $G$ is torsion free, $\ZZZ e \subset G$ and that $G/\ZZZ e$ is a
torsion group. Under these conditions it was proved in \cite{WinterToms:ssa} that
$G$ is isomorphic to either $\ZZZ$ or to a subgroup of $H$ of $\mathbb{Q}$ with
the property that $1/n^2 \in H$ whenever $1/n\in H$ for some nonzero integer $n$.
Indeed, for every  $x\in G$, $x \neq 0$, there is a unique pair of relatively
prime integers $m$ and $n$ with $n>0$ such that $nx=me$. One verifies immediately
that $\gamma:G \to \mathbb{Q}$, $\gamma(x)=m/n$, if $x\neq 0$ and $\gamma(0)=0$
is an injective morphism of groups. Moreover, if $x\in G$ satisfies
$\gamma(x)=1/n$, then $y=\theta^{-1}(x\ot x)$ satisfies $\gamma(y)=1/n^2$.
\end{proof}
\emph{Proof of Theorem~\ref{Theorem:ssa-Kirchberg-algebras}}
\begin{proof} Let $D$ be isomorphic to either $\OOO_2$,  $\OOO_\infty$, or
$\OOO_\infty\ot U$ where $U$ is a unital UHF algebra of infinite type. Let $X$ be
a finite dimensional compact metrizable space and let $A$ be a unital separable continuous
$C(X)$-algebra with all fibers isomorphic to $D$.
 Let $\varphi: D \to A\otimes D$ be defined
by $\varphi( d)=1_A\otimes d$ for all  $d\in D$ and let
$\widetilde{\varphi}:C(X)\otimes D \to A\otimes D$ be the $C(X)$-linear extension
of $\varphi$. By the UCT, every unital $*$-endomorphism of $D$ is a
KK-equivalence. Therefore each $\varphi_x:D \to A(x)\otimes D\cong D$ is a
KK-equivalence. It follows that
 $\KKX(\widetilde{\varphi})$ satisfies the assumptions of
Theorem~\ref{Theorem:general-KK-for Kirchberg algebras} and hence $A\otimes D$ is
isomorphic to $C(X)\otimes D$. We conclude the first half of the proof by invoking a recent result
of Hirshberg, R{\o}rdam and Winter
\cite[Thm.~4.3]{HirshbergRordamWinter:absorb-ssa} which shows that if $D$ is
strongly selfabsorbing, then $A\otimes D\cong A$.

Conversely, suppose that $D$ is a unital Kirchberg algebra which has the
automatic triviality property and satisfies the UCT. Let $Y$ be a  finite dimensional compact metric
space and let $SY$ be its unreduced suspension. The locally trivial
$C(SY)$-algebras with fibers $D$ and structure group $\Aut(D)^0$ are classified
by the homotopy classes $[Y,\Aut(D)^0]$ (see \cite{Hus:fibre}) and hence
$[Y,\Aut(D)^0]$ must reduce to a singleton since $D$ has the automatic triviality
property. We conclude the proof by applying Proposition~\ref{Tor}.
\end{proof}
\section{$C(X)$-algebras and the Universal Coefficient
Theorem}\label{sec-uct}

Let us recall the notion of \emph{category} of a  $C(X)$-algebra with respect
   to a
    class $\mathcal{C}$ of C*-algebras
    \cite{Dad:bundles-fdspaces}.
    A  $C(Z)$-algebra $E$ satisfies $\mathrm{cat}_\mathcal{C}(E)=0$
  if there is a
finite partition of $Z$
    into   closed subsets $Z_1,\dots,Z_r$ ($r\geq
    1$) and there exist C*-algebras $D_1,\dots,D_r$ in $ \mathcal{C}$
     such that $E\cong\bigoplus_{i=1}^r C(Z_i)\ot D_i$. We write
    $\mathrm{cat}_\mathcal{C}(A)\leq n$ if there are closed
    subsets $Y$ and $Z$ of $X$  with
    $X=Y\cup Z$ and there exist  a  $C(Y)$-algebra $B$ with
     $\mathrm{cat}_\mathcal{C}(B)\leq n-1$ and
    a
     $C(Z)$-algebra $E$ with $\mathrm{cat}_\mathcal{C}(E)=0$
   and  a  $*$-monomorphism of $C(Y\cap Z)$-algebras
    $\gamma:E(Y\cap Z)\to B(Y\cap
Z)$ such that $A$ is isomorphic to
$$B\oplus_{\,\pi,\gamma\pi_{}} E=\{(b,e)\in B \oplus E:\,
\pi^Y_{Y\cap Z}(b)=\gamma\pi^Z_{Y\cap Z}(e)\}.$$ By definition
$\mathrm{cat}_{\mathcal{C}}(A)=n$ if $n$ is the smallest number with the property
that $\mathrm{cat}_{\mathcal{C}}(A)\leq n.$ If  no such $n$ exists, then
$\mathrm{cat}_{\mathcal{C}}(A)=\infty.$ One has an exact sequence
 \begin{equation}\label{1}
\xymatrix{
 {0}\ar[r]& {\{b\in B:\pi_{Y\cap Z}(b)=0\}}\ar[r]
                & B\oplus_{\,\pi,\gamma\pi_{}} E \ar[r]^{\,\,\quad\pi_Z}&{E}\ar[r]&0
}
\end{equation}
\begin{lemma}\label{uct-help}
Let $A$ be a
    $C(X)$-algebra such that
$\mathrm{cat}_\mathcal{C}(A)=n <\infty$ where $\mathcal{C}$ is the class of all
 Kirchberg algebras satisfying
the UCT. Then $A$ satisfies the UCT.
\end{lemma}
\begin{proof} We shall
 prove
 by induction on $n$ that if $\mathrm{cat}_{\mathcal{C}}(A)\leq n$,
 then  $A$ and all its  closed two-sided ideals satisfy the UCT. If $n=0$, then
  $A\cong\oplus_{i} C(Z_i)\ot D_i$ and all its closed two-sided ideals satisfy the UCT since
 each $D_i$ is simple and satisfies the UCT. By a result of \cite{RosSho:UCT}, if
 two out of  three  separable nuclear C*-algebras
 in a short exact sequence
 satisfy the UCT, then all three of them satisfy the UCT.
 For the inductive step we use the
 exact sequence~\eqref{1}, with $E$
  elementary and $\mathrm{cat}_{\mathcal{C}}(B)\leq n-1$.
\end{proof}
\begin{theorem}[\cite{Dad:uct}]\label{mainn}
Let $A$ be a nuclear separable C*-algebra. Assume that  for any finite set
$\fset \subset A$ and any $\ep>0$ there is a C*-subalgebra $B$ of $A$
satisfying the UCT and such that $\fset\subset_\ep B$. Then $A$ satisfies the
UCT.
\end{theorem}
\begin{proof} For the convenience of the reader  we sketch an  alternate proof in the case when $B$
is nuclear. It is just this case that is needed in the sequel. By assumption,
$A$ admits  an exhaustive sequence $(A_n)$ consisting of nuclear separable
C*-subalgebras which satisfy the UCT. We may assume that $A$ is unital and its
unit is contained in each $A_n$. Let us replace the pair $A_n\subseteq A$ by
$A_n\otimes p\OOO_\infty p\subseteq A\otimes p\OOO_\infty p$ with $p$ as in
Lemma~\ref{O2}. If we use the map $\theta:A\otimes p\OOO_\infty
p\hookrightarrow \OOO_2\subset 1_A\ot p\OOO_\infty p$ and $s_1,s_2\in 1_A\ot
p\OOO_\infty p$ to construct $\varphi:A\otimes p\OOO_\infty p \to A\otimes
p\OOO_\infty p$ as in Theorem~\ref{C(X)-K-version}, then $\varphi(B\otimes
p\OOO_\infty p)\subset B\otimes p\OOO_\infty p$ for any subalgebra $B$ of $A$.
Therefore we can make
 the construction $A\mapsto A^\sharp$  functorial with
respect to subalgebras. This shows that  $A^\sharp$ admits  an exhaustive
sequence $(A_n^\sharp)$ consisting of nuclear separable C*-subalgebras which
satisfy the UCT since each $A_n^\sharp$ is KK-equivalent to $A_n$. We can write
each $A_n^\sharp$ as an inductive limit of a sequence of Kirchberg algebras
satisfying the UCT and having finitely generated K-theory groups \cite{Ror:ency}. Those algebras
are weakly semiprojective (\cite{Lin:wsp},\cite{Spi:semiproj}; see also   \cite[Thm.~3.11]{Dad:bundles-fdspaces} for a short proof). Thus $A^\sharp$
admits
 an exhaustive sequence $(B_n)$  consisting of weakly
semiprojective
 C*-algebras which satisfy the UCT.
By a standard
 perturbation argument (\cite{Lor:lifting}) we see
 that
  $A^\sharp$ is isomorphic to the inductive limit of a subsequence
$(B_{i_n})$ of $(B_n)$ and hence  $A^\sharp$ satisfies the UCT
\cite{RosSho:UCT}. Therefore $A$ satisfies the UCT since it is KK-equivalent to
$A^\sharp$ .
\end{proof}
\emph{Proof of Theorem~\ref{intro-uct}}
\begin{proof}
 Let $A$ be as in the statement and
consider the open set $Y=\{x\in X: A(x)\neq 0\}$.  By replacing $X$ by $Y$ and
viewing $A\cong C_0(Y)A$ as a $C(Y)$-algebra we may assume that all the fibers of
$A$ are nonzero. Let $X^+$ be the one-point compactification of $X$. Then
$C(X^+)$ is separable by \cite[Lemma~2.2]{Dad:bundles-fdspaces}. By
{\cite[Prop.~3.2]{Blanchard:Hopf}}, there is a unital $C(X^+)$-algebra $A^+$
which contains $A$ as an ideal and such that $A^+/A\cong C(X^+)$. Thus we have
reduced the proof to the case when  $X$ is compact and metrizable and $A$ is
unital.
 By Theorem~\ref{C(X)-K-version} we may assume that the
fibers of $A$
are Kirchberg
     C*-algebras satisfying the UCT. By \cite[Thm.~4.6]{Dad:bundles-fdspaces},
     $A$ admits an exhaustive sequence $(A_k)$ such that each $A_k$ verifies
     the assumptions of Lemma~\ref{uct-help} and hence
    $A_k$ satisfies the UCT.
     We conclude the proof by applying Theorem~\ref{mainn}.
Let us note that the above proof only requires a weaker version of
Theorem~\ref{C(X)-K-version} which states that $\Phi:A \to A^\sharp$  and each
$\Phi_x$ are KK-equivalences. Its proof requires only the usual
$\varprojlim^1$-sequence for $KK$-theory.
\end{proof}

\bibliographystyle{abbrv}
\small{\bibliography{operator}}
\end{document}